%
%
\input amstex.tex
\documentstyle{amsppt}
\magnification=1200
\hsize=6.5truein
\baselineskip=13pt

\def\R{{\bold R}}
\def\N{{\bold N}}
\def\C{{\bold C}}
\def\Z{{\bold Z}}
\def\Zp{{\bold Z}_+}

\topmatter
\title A $q$-analogue of Graf's addition formula \\
for the Hahn-Exton $q$-Bessel function\endtitle
\rightheadtext{Graf's addition formula for Hahn-Exton $q$-Bessel function}
\author H.T. Koelink and R.F. Swarttouw\endauthor
\affil  Katholieke Universiteit Leuven
and Vrije Universiteit Amsterdam\endaffil
\address Departement Wiskunde, Katholieke Universiteit Leuven,
Celestijnenlaan 200 B, B-3001 Heverlee, Belgium\endaddress
\email erik\%twi\%wis\@cc3.KULeuven.ac.be\endemail
\address Faculteit Wiskunde en Informatica, Vrije Universiteit,
De Boelelaan 1081, 1081 HV Amsterdam, the Netherlands\endaddress
\email rene\@cs.vu.nl\endemail
\thanks HTK is supported by a NATO-Science Fellowship of the
Netherlands Organization for Scientific Research (NWO). \endthanks
\abstract An addition and product formula for the Hahn-Exton $q$-Bessel
function, previously obtained by use of a quantum group
theoretic interpretation, are proved analytically. A (formal)
limit transition to the Graf addition formula and corresponding
product formula for the Bessel function is given.
\endabstract
\endtopmatter
\document

\head 1. Introduction\endhead

A classical result for the Bessel function
$$
J_\nu(z) = \sum_{k=0}^\infty {{(-1)^k(z/2)^{\nu+2k}}\over{k!\,
\Gamma(\nu+k+1)}}
\tag1.1
$$
is the addition formula
$$
J_\nu\Bigl( \sqrt{x^2+y^2-2xy\cos \psi}\Bigr)
\left( {{x-ye^{-i\psi}}\over{x-ye^{i\psi}}} \right)^{\nu/2} =
\sum_{m=-\infty}^\infty J_{\nu+m}(x)J_m(y)e^{im\psi},
\tag1.2
$$
$\vert ye^{\pm i\psi} \vert < \vert x\vert$,
due to Graf (1893), cf. \cite{19, \S 11.3(1)}, for general $\nu$, and
due to Neumann (1867) for $\nu=0$, cf. \cite{19, \S 11.2(1)}.
In case $\nu\in\Z$ the conditions on $x$, $y$ in \thetag{1.2}
can be removed.

Several $q$-analogues of the Bessel function \thetag{1.1} have been
studied. The oldest $q$-analogues have been introduced by Jackson in
1903-1905, cf. the references in Ismail \cite{4}. The $q$-Bessel
function studied in this note is the so-called Hahn-Exton $q$-Bessel
function which has been introduced by Hahn (1949) for
a special case and by Exton (1978) in full generality,
cf. references in Koornwinder and Swarttouw \cite{13}.

The Hahn-Exton $q$-Bessel function is defined by
$$
J_\alpha (z;q) = z^\alpha {{(q^{\alpha+1};q)_\infty}\over{(q;q)_\infty}}
\, {}_1 \varphi_1 \left( {{0}\atop{q^{\alpha+1}}};q,qz^2 \right).
\tag1.3
$$
Here $q\in (0,1)$, $(a;q)_0=1$, $(a;q)_k=\prod_{i=0}^{k-1}(1-aq^i)$ for
$k\in\N$, $(a;q)_\infty=\lim_{k\to\infty}(a;q)_k$ and the $q$-hypergeometric
function is defined by
$$
{}_r\varphi_s \left( {{a_1,a_2,\ldots,a_r}\atop{b_1,\ldots,b_s}};q,z\right) =
\sum_{k=0}^\infty {{(a_1;q)_k(a_2;q)_k\ldots (a_r;q)_k}\over
{(q;q)_k(b_1;q)_k\ldots (b_s;q)_k}} z^k
\Bigl( (-1)^kq^{{1\over 2}k(k-1)}\Bigr)^{(s-r+1)}.
$$
The notation for $q$-shifted factorials and $q$-hypergeometric series is
taken from the book by Gasper and Rahman \cite{3} to which the reader
is referred for more information on this subject.

The goal of this note is to prove the formula
$$
\eqalign{
&\qquad \qquad J_\nu(Rq^{{1\over 2}(y+z+\nu)};q)
J_{x-\nu}(q^{{1\over 2}z};q)\cr
=& \sum_{k=-\infty}^\infty J_k(Rq^{{1\over 2}(x+y+k)};q)
J_{\nu+k}(Rq^{{1\over 2}(y+k+\nu)};q)J_x(q^{{1\over 2}(z-k)};q).\cr}
\tag1.4
$$
This formula is valid for $z\in\Z$, $R,x,y,\nu\in\C$ satisfying
$q^{1+\Re (x)+\Re (y)}\vert R\vert^2 < 1$, $\Re (x) > -1$,
and $R \not= 0$.

The formula \thetag{1.4} has originally been derived
for $\nu,x,y\in\Z$, $R>0$
by Koelink using the interpretation of the Hahn-Exton $q$-Bessel function
as matrix elements of irreducible unitary representations of the
quantum group of plane motions, cf. \cite{7, \S 3.5}. This
quantum group theoretic interpretation of the Hahn-Exton $q$-Bessel
function is due to Vaksman and Korogodski\u\i\ \cite{16}.
The quantum group theoretic derivation of \thetag{1.4}
in \cite{7}
is modelled on the group theoretic derivation of Graf's addition
formula \thetag{1.2} as presented by
Vilenkin and Klimyk \cite{18,
\S 4.1.4(2)}, so we call \thetag{1.4} a
$q$-analogue of Graf's addition formula for the Hahn-Exton $q$-Bessel
function.

In section 4 we present a (formal) limit transition of
\thetag{1.4} to \thetag{1.2} as
$q$ tends to $1$. The approach employed is based on a theorem
presented by Van Assche and Koornwinder \cite{17} which has been
used to show that the addition formula for the little
$q$-Legendre polynomial \cite{12} tends to the addition formula for the
Legendre polynomial. The case $\nu=0$, $R=q^{-{1\over 2}}$ of the addition
formula \thetag{1.4} can be obtained by taking a (formal)
limit in Koornwinder's \cite{12} addition formula for the little
$q$-Legendre polynomials using the limit transition of the little
$q$-Jacobi polynomials to the Hahn-Exton $q$-Bessel function, cf.
\cite{13, prop. A.1}.

Kalnins, Miller and Mukherjee \cite{6}
have given another derivation of the
addition formula \thetag{1.4}.
They consider representations of
the Lie algebra of the group of orientation and distance preserving motions
of the plane. Instead of exponentiating the representations using the
exponential function, they use a $q$-analogue of the exponential function. In
a particular case the matrix elements can be expressed in terms of the
Hahn-Exton $q$-Bessel functions. They give a decomposition of the tensor
product (not the standard tensor product, but one closely related to
quantum groups) of two representations and the corresponding Clebsch-Gordan
coefficients are in terms of the Hahn-Exton $q$-Bessel functions. Comparison
of matrix coefficients yields the addition formula
\thetag{1.4}.

More results on Graf type addition formulas for the Hahn-Exton $q$-Bessel
function, the Jackson $q$-Bessel function and other $q$-analogues of
the Bessel function can be found in \cite{1, 2,
5, 6, 7, 8, 9, 13}.

We start the proof of the $q$-analogue of Graf's addition formula
\thetag{1.4}
for the Hahn-Exton $q$-Bessel function by proving the product formula
$$
\eqalign{
&\qquad (-1)^m q^{-{1\over 2}m} J_m(Rq^{{1\over 2}(x+y)};q)
J_{\nu-m}(Rq^{{1\over 2}(y+\nu-m)};q)\cr
&= \sum_{z=-\infty}^\infty q^z J_x(q^{{1\over 2}(m+z)};q)
J_{x-\nu}(q^{{1\over 2}z};q) J_\nu(Rq^{{1\over 2}(y+\nu+z)};q) \cr}
\tag1.5
$$
valid for $m\in\Z$, $R,x,y,\nu\in\C$ satisfying
$\Re (x) > -1$, $q^{1+\Re (x)+\Re (y)}\vert R\vert^2 < 1$
and $R \not= 0$. This is done in \S 2. The proof of the addition
formula \thetag{1.4} is a direct consequence of
\thetag{1.5} as is shown in \S 3.
The product formula in the case $\nu=m=0$, $R=q^{-{1\over 2}}$ is
mentioned (without proof) by Vaksman and Korogodski\u\i\
\cite{16, p.177}.

In \S 5 we will show that \thetag{1.5} is a $q$-analogue
of the product formula for Bessel functions,
$$
J_{\nu+m}(x)J_m(y)=
{{1}\over{2\pi}} \int_0^{2\pi}
J_\nu\Bigl( \sqrt{x^2+y^2-2xy\cos \psi}\Bigr)
\left( {{x-ye^{-i\psi}}\over{x-ye^{i\psi}}} \right)^{\nu/2}
e^{-im\psi}\,d\psi .
\tag1.6
$$
The product formula is a direct consequence of Graf's addition formula
\thetag{1.2}.

Acknowledgement. We thank the referee, A.P.~Magnus and W.~Van~Assche
for constructive remarks concerning \S 4.

\head 2. Proof of the product formula
\endhead

In this section we present an analytic proof of the product
formula \thetag{1.5}. The proof uses two known formulas for
the Hahn-Exton $q$-Bessel function previously obtained
by Koornwinder and Swarttouw \cite{13} and Swarttouw \cite{14}.

The proof starts with the following formula, valid for
$\vert sxy \vert < 1$ and $m\in\Z$;
$$
\eqalign{
&\qquad y^m  {{(s^{-1}xy^{-1};q)_\infty
(q^{m+1};q)_\infty}\over{(q^ms^{-1}xy^{-1};q)_\infty (q;q)_\infty}}
\, {}_2\varphi_1 \left( {{q^ms^{-1}xy^{-1},s^{-1}yx^{-1}}\atop{
q^{m+1}}};q,sxy\right) \cr
=  & \sum_{z=-\infty}^\infty s^z y^{m+z} {{(y^2;q)_\infty}
\over{(q;q)_\infty}} \,{}_1 \varphi_1 \left( {{0}\atop{y^2}};
q,q^{m+z+1}\right) x^z {{(x^2;q)_\infty}\over{(q;q)_\infty}}
\,{}_1 \varphi_1 \left( {{0}\atop{x^2}};q,q^{z+1}\right), \cr}
\tag2.1
$$
cf. \cite{13, (4.5), (2.3)}.

In \thetag{2.1} we take $x=q^{{1\over 2}(x-\nu+1)}$,
$y=q^{{1\over 2}(x+1)}$ and $s=q^{{1\over 2}\nu+n}$ for $n\in\Zp$
to get
$$
\eqalign{
q^{{1\over 2}mx} & {{(q^{-\nu-n};q)_\infty
(q^{m+1};q)_\infty}\over{(q^{m-\nu-n};q)_\infty (q;q)_\infty}}
\, {}_2\varphi_1 \left( {{q^{m-\nu-n},q^{-n}}\atop{
q^{m+1}}};q,q^{x+n+1}\right) \cr
=  & \sum_{z=-\infty}^\infty q^{z(n+1+{1\over 2}\nu)}
J_x(q^{{1\over 2}(m+z)};q) J_{x-\nu}(q^{{1\over 2}z};q)\cr}
\tag2.2
$$
valid for $m\in\Z$, $\nu\in\C$, $\Re (x) > -1$ and all $n\in\Zp$.
In \thetag{2.2} we use the notation \thetag{1.3}.

Multiply both sides of \thetag{2.2} by
$$
{{(q^{\nu+1};q)_\infty}\over{(q;q)_\infty}} q^{{1\over 2}\nu(\nu+y)}
{{(-1)^nq^{{1\over 2}n(n+1)}}\over{(q^{\nu+1};q)_n(q;q)_n}} q^{n(y+\nu)}
R^{\nu+2n}
$$
and sum from $n=0$ to $\infty$. After interchanging the summations
over $z\in\Z$ and $n\in\Zp$, which is justified
for $\Re (x) > -1$ and $q^{1+\Re (x) + \Re (y)} \vert R\vert^2 < 1$,
cf. proposition A.1, we find
$$
\eqalign{
& \qquad\qquad \sum_{z=-\infty}^\infty q^z J_x(q^{{1\over 2}(m+z)};q)
J_{x-\nu}(q^{{1\over 2}z};q) J_\nu(Rq^{{1\over 2}(y+\nu+z)};q)\cr
& \qquad\qquad
= q^{{1\over 2}mx}{{(q^{m+1};q)_\infty(q^{\nu+1};q)_\infty}\over
{(q;q)_\infty(q;q)_\infty}} q^{{1\over 2}\nu(\nu + y)}R^\nu \cr
&\times \sum_{n=0}^\infty {{(q^{-\nu-n};q)_\infty}\over
{(q^{m-\nu-n};q)_\infty}}
{{(-1)^nq^{{1\over 2}n(n+1)}}\over{(q^{\nu+1};q)_n(q;q)_n}} q^{n(y+\nu)}
R^{2n} \, {}_2\varphi_1 \left( {{q^{m-\nu-n},q^{-n}}\atop{
q^{m+1}}};q,q^{x+n+1}\right). \cr}
\tag2.3
$$
We use for $m\in\Z$, $n\in\Zp$,
$$
{{(q^{\nu+1};q)_\infty (q^{-\nu-n};q)_\infty}\over
{(q^{\nu+1};q)_n (q^{m-\nu-n};q)_\infty}} =
(-1)^m q^{{1\over 2}m(m-1)-m(\nu+n)}
{{(q^{\nu-m+1};q)_\infty}\over{(q^{\nu-m+1};q)_n}}
$$
to see that the right hand side of \thetag{2.3} can be
rewritten as
$$
\eqalign{
& (-1)^m q^{-{1\over 2}m}
{{(q^{m+1};q)_\infty(q^{\nu-m+1};q)_\infty}\over
{(q;q)_\infty(q;q)_\infty}}
q^{{1\over 2}(\nu-m)(\nu-m + y)+ {1\over 2}m(x+y)}R^\nu \cr
&\times \sum_{n=0}^\infty
{{(-1)^nq^{{1\over 2}n(n+1)}}\over{(q^{\nu-m+1};q)_n(q;q)_n}} q^{n(y+\nu-m)}
R^{2n} \, {}_2\varphi_1 \left( {{q^{m-\nu-n},q^{-n}}\atop{
q^{m+1}}};q,q^{x+n+1}\right). \cr}
\tag2.4
$$

Now we apply the product formula for the Hahn-Exton $q$-Bessel function,
cf. \cite{14, (3.1)}, which is valid
for $a,b,x,\mu,\nu\in\C$ provided $abx\neq 0$. Explicitly,
$$
\eqalign{
&J_\nu(ax;q)J_\mu(bx;q)= {{(q^{\nu+1};q)_\infty(q^{\mu+1};q)_\infty}\over
{(q;q)_\infty(q;q)_\infty}}a^\nu b^\mu x^{\nu+\mu} \cr
&\qquad\times \sum_{n=0}^\infty {{(-1)^n(bx)^{2n}q^{{1\over 2}n(n+1)}}\over
{(q^{\mu+1};q)_n(q;q)_n}}\, {}_2 \varphi_1 \left(
{{q^{-n},q^{-n-\mu}}\atop{q^{\nu+1}}};q,q^{\mu+n+1}{{a^2}\over{b^2}}\right).
\cr}
\tag2.5
$$
From \thetag{2.5} we see that \thetag{2.4} equals
$$
(-1)^m q^{-{1\over 2}m} J_m(Rq^{{1\over 2}(x+y)};q)
J_{\nu-m}(Rq^{{1\over 2}(y+\nu-m)};q),
$$
which proves the product formula \thetag{1.5}.

\head 3. Proof of the addition formula
\endhead

The proof of the addition formula \thetag{1.4} uses the
orthogonality relations
$$
\sum_{m=-\infty}^\infty q^{m+z} J_x(q^{{1\over 2}(z+m)};q)
J_x(q^{{1\over 2}(l+m)};q) = \delta_{z,l},
\tag3.1
$$
for $z,l\in\Z$, $\Re (x) > -1$, cf. \cite{13, (2.11)}. In
\cite{13, \S 3} it is shown that \thetag{3.1} can be
viewed as a $q$-analogue of the Hankel transform (or the Fourier-Bessel
integral).

Multiply both sides of the product formula \thetag{1.5}
by $q^mJ_x(q^{{1\over 2}(l+m)};q)$ for $l\in\Z$ and sum over
$m$ from $-\infty$ to $\infty$. We interchange the summations over
$m\in\Z$ and $z\in\Z$, which is allowed for
$\Re (x) > -1$, $q^{1+\Re (x)+\Re (y)}\vert R\vert^2 < 1$,
cf. proposition A.2. An application of the orthogonality relations
\thetag{3.1} and replacing $m$ by $-m$ yields
$$
\eqalign{
&\qquad \qquad J_\nu(Rq^{{1\over 2}(y+l+\nu)};q)
J_{x-\nu}(q^{{1\over 2}l};q)\cr
=& \sum_{m=-\infty}^\infty (-1)^m q^{-{1\over 2}m}
J_{-m}(Rq^{{1\over 2}(x+y)};q)
J_{\nu+m}(Rq^{{1\over 2}(y+m+\nu)};q)J_x(q^{{1\over 2}(l-m)};q).\cr}
$$
Since $J_{-n}(z;q)=(-1)^nq^{{1\over 2}n}J_n(zq^{{1\over 2}n};q)$,
$n\in\Z$, cf. \cite{13, (2.6)}, the addition formula
\thetag{1.4} is proved.

\head 4. The limit case $q\uparrow 1$ of the addition formula
\endhead

In this section we present a (formal) limit transition of the
$q$-analogue of the addition formula \thetag{1.4}
for the Hahn-Exton $q$-Bessel function to Graf's addition
formula \thetag{1.2} for the Bessel function.

First we recall the $q$-gamma function, cf. \cite{3, \S 1.10},
$$
\Gamma_q(z) = {{(q;q)_\infty}\over{(q^z;q)_\infty}}(1-q)^{1-z},
\qquad \lim_{q\uparrow 1}\ \Gamma_q(z) = \Gamma(z).
\tag4.1
$$
The $q$-gamma function can be used to see that (formally)
$$
\lim_{q\uparrow 1} J_\nu (z(1-q)/2;q) = J_\nu(z).
\tag4.2
$$
See \cite{13, App. A} for a rigorous limit result of this type.

In order to obtain Graf's addition formula \thetag{1.2}
from \thetag{1.4} as $q\uparrow 1$ we have to consider
the quotient of two Hahn-Exton $q$-Bessel functions, which we rewrite
as an infinite sum of quotients of two Hahn-Exton $q$-Bessel functions
of equal order. Explicitly,
$$
{{J_x(q^{{1\over 2}(z-k)};q)}\over
{J_{x-\nu}(q^{{1\over 2}z};q)}} = q^{{1\over 2}\nu(z-k)}
\sum_{m=0}^\infty {{(q^\nu ;q)_m}\over{(q;q)_m}}
q^{m(1+{1\over 2}(x-\nu))}
{{J_{x-\nu}(q^{{1\over 2}(z-k+m)};q)}\over
{J_{x-\nu}(q^{{1\over 2}z};q)}}.
\tag4.3
$$
Note that the Hahn-Exton $q$-Bessel functions on the
right hand side of \thetag{4.3} only differ in the argument
by a (half-)integer power of $q$. Equation \thetag{4.3}
holds for $J_{x-\nu}(q^{{1\over 2}z};q)\not= 0$ and
$\Re (x-\nu)>-1$. It can be proved by substituting the series
representation \thetag{1.3} for the Hahn-Exton $q$-Bessel
function $J_{x-\nu}$ in the nominator on the right hand side of
\thetag{4.3}, interchanging summations and using the
$q$-binomial formula.

The limit $q\uparrow 1$ in \thetag{4.3}
can be taken with the help of the
following proposition. The proof of the proposition is
an easy adaptation of the proof of \cite{17, thm. 1},
which we will not give.

\proclaim{Proposition 4.1}
Suppose $\{ p_k(x;n)\mid k\in\Zp, n\in\N\}$ is a series of
functions satisfying a recurrence relation
$$
x^2p_k(x;n) = a_{k+1,n}p_{k+1}(x;n) + b_{k,n}p_k(x;n)
+ a_{k,n}p_{k-1}(x;n)
\tag4.4
$$
with $a_{k,n}>0$, $b_{k,n}\in\R$ and initial conditions
$p_0(x;n)=f(x;n)$, $p_1(x;n)=g(x;n)$. Assume that
the zeros of $p_k(x;n)$ are real for all $k\in\Zp$ and all
$n\in\N$ and that
$$
\Bigl\vert {{p_{k-1}(x;n)}\over{p_k(x;n)}}\Bigr\vert <
{{C}\over{\delta}}, \qquad\quad \forall\, k,n\in\N,
\tag4.5
$$
uniformly for $x\in K$ for any compact
subset $K$ of $\C\backslash\R$ with $d(K,\R) > \delta$.
Moreover, assume that
$$
\eqalign{
\lim_{n\to\infty} a_{n,n} = A > 0, \qquad
&\lim_{n\to\infty} b_{n,n} = B \in\R, \cr
\lim_{n\to\infty} a_{k,n}^2 - a_{k-1,n}^2 = 0, \qquad
&\lim_{n\to\infty} b_{k,n} - b_{k-1,n} = 0,\cr}
$$
uniformly in $k$.

Then we have
$$
\lim_{n\to\infty} {{p_{n+1}(x;n)}\over{p_n(x;n)}} =
\rho\left( {{x^2-B}\over{2A}}\right)
\tag4.6
$$
uniformly on compact subsets of $\C\backslash \R$ with
$\rho(x)=x+\sqrt{x^2-1}$ (and the square root is defined by
$\vert \rho(x)\vert >1$ for $x\notin [-1,1]$).
\endproclaim

Next we apply proposition 4.1 to
$$
p_k(x;n)= (-1)^k q^{-{1\over 2}k} J_{2n\alpha+\beta}(xq^{-{1\over 2}k};q)
\tag4.7
$$
with $q$ replaced by $c^{1/n}$ for a fixed $c\in (0,1)$. We assume
$\alpha>0$ and $\beta >-1$, so that all zeros of $p_k(x;n)$ are real
by \cite{10, cor. 3.2}.
It follows
from the second order $q$-difference equation for the Hahn-Exton
$q$-Bessel function, cf. \cite{15, (17), (18)}, that
\thetag{4.4} is satisfied with
$$
a_{k,n}=q^{k-{1\over 2}+n\alpha+{1\over 2}\beta}, \qquad
b_{k,n}=q^k(1+q^{2n\alpha+\beta})
$$
and $q=c^{1/n}$. Then the conditions on the coefficients $a_{k,n}$
and $b_{k,n}$ of proposition 4.1 are easily verified. Moreover,
$A=c^{1+\alpha}>0$ and $B=c(1+c^{2\alpha})\in\R$. It remains
to prove that \thetag{4.5} holds.

\proclaim{Lemma 4.2}
With $p_k(x;n)$ defined by {\rm \thetag{4.7}} for $\alpha>0$
and $\beta >-1$, the estimate {\rm \thetag{4.5}} holds.
\endproclaim

\demo{Proof} Consider the Wall polynomials for $0<b<1$,
$$
w_p(x;b;q) = (-1)^p \sqrt{ {{(b;q)_p}\over{b^p(q;q)_p}} }
\, {}_2\varphi_1 \left(  {{q^{-p},0}\atop{b}};q,x\right).
$$
From \cite{17, (2.5), cor. 1} we get for $p\in\N$
$$
\Bigl\vert {{w_{p-1}(x;b;q)}\over{w_p(x;b;q)}}\Bigr\vert \leq
{{q^p\sqrt{b(1-q^p)(1-bq^{p-1})}}\over{\delta'}} \leq
{{q^p\sqrt{b}}\over{\delta'}}
$$
for all $x\in K'$, $K'$ compact in $\C \backslash [0,1]$,
$\delta' = d(K',[0,1])$. Replace $b$, $x$ by $q^{\mu+1}$,
$y^2q^{p-m}$ and take the limit $p\to\infty$. The Wall
polynomials tend to the Hahn-Exton $q$-Bessel function
uniformly, cf. \cite{13, prop. A.1}, and we
get
$$
\Bigl\vert {{J_{\mu}(yq^{-{1\over 2}m};q)}\over
{J_\mu(yq^{-{1\over 2}(m+1)};q)}}\Bigr\vert \leq
{{q^{m+{1\over 2}\mu}}\over{\delta}}, \qquad \mu>-1
\tag4.8
$$
for all $y\in K$, $K$ compact set of $\C \backslash \R$ with
$\delta=d(K,\R)$. The required estimate \thetag{4.5}
for $p_k(x;n)$ defined by \thetag{4.7} follows
from  \thetag{4.8}. \qed
\enddemo

Consequently, by proposition 4.1,
$$
\lim_{n\to\infty}
{{J_{2n\alpha+\beta}(xc^{-{1\over{2n}}}c^{-{1\over 2}};c^{1\over n})}\over
{J_{2n\alpha+\beta}(xc^{-{1\over 2}};c^{1\over n})}} =
\rho = \rho(x;\alpha,c), \qquad
\rho +{1\over \rho} = {{c(1+c^{2\alpha})-x^2}\over{c^{1+\alpha}}}.
\tag4.9
$$
In \thetag{4.3} we replace $x$, $z$, $q$ by $2n\alpha$, $2n\gamma$,
$c^{1\over n}$, and we take $n\to\infty$. Iterating \thetag{4.9}
and interchanging limit and summation shows that
\thetag{4.3} tends to
$$
c^{\nu\gamma} \sum_{m=0}^\infty {{(\nu)_m}\over{m!}} c^{m\alpha}\rho^{k-m} =
c^{\nu\gamma}\rho^k (1-c^\alpha \rho^{-1})^{-\nu}
\tag4.10
$$
with $\rho+1/\rho = (1+c^{2\alpha}-c^{2\gamma})/c^{\alpha}$, $\nu<1$ and
$c^\gamma\in\C\backslash\R$. Interchanging limit and summation can be
justified using the estimate \thetag{4.8},
\cite{11, lemma A.1} and dominated convergence.

Finally, divide both sides of the addition formula
\thetag{1.4} by $J_{x-\nu}(q^{{1\over 2}z};q)$.
Replace $x$, $z$, $q$ by $2n\alpha$, $2n\gamma$,
$c^{1\over n}$ as before and $y$, $R$ by $2n\eta$, $R(1-q)/2$.
Take $n\to\infty$ and use \thetag{4.2},
\thetag{4.10} to see that \thetag{1.4}
formally tends to
$$
\sum_{k=-\infty}^\infty \rho^k J_k(Rc^{\alpha+\eta})J_{\nu+k}(Rc^\eta) =
c^{-\nu\gamma}(1-\rho^{-1}c^\alpha)^\nu J_\nu(Rc^{\eta+\gamma}),
$$
which is easily rewritten as Graf's addition formula
\thetag{1.2} for the Bessel function.

\head 5. The limit case $q\uparrow 1$ of the product formula
\endhead

A (formal) limit transition of the product formula
\thetag{1.5} for the Hahn-Exton $q$-Bessel
function to product formula \thetag{1.6} for the Bessel
function is presented in this section.

We first rewrite the product formula \thetag{1.5}.
Use $J_m(z;q)=(-1)^mq^{-{1\over 2}m} J_{-m}(zq^{-{1\over 2}m};q)$,
replace $m$ by $-m$, and use the formula
$$
J_{x-\nu}(z;q) = z^{-\nu}\sum_{k=0}^\infty {{(q^{-\nu};q)_k}\over{(q;q)_k}}
q^{k(1+{1\over 2}x)} \, J_x(zq^{{1\over 2}k};q)
$$
for $\Re (x)>-1$, and the notation $\int_0^\infty f(z)dm_q(z)=
\sum_{z=-\infty}^\infty q^zf(q^{{1\over 2}z})$
to write \thetag{1.5} as
$$
\eqalign{
&\qquad\qquad
q^m J_{\nu+m}(Rq^{{1\over 2}(y+\nu+m)};q)J_m(Rq^{{1\over 2}(x+y+m)};q) = \cr
&\sum_{k=0}^\infty {{(q^{-\nu};q)_k}\over{(q;q)_k}}
q^{k(1+{1\over 2}x)} \int_0^\infty z^{-\nu}J_{\nu}(zRq^{{1\over 2}(y+\nu)};q)
J_x(zq^{-{1\over 2}m};q)J_x(zq^{{1\over 2}k};q) dm_q(z). \cr}
\tag5.1
$$

In order to calculate the limit of the integral on the right hand side of
\thetag{5.1} we assume $x=2n\alpha$, $q=c^{1\over n}$, for
$\alpha > 0$, $c\in (0,1)$, so that we can use \thetag{4.7}
(with $\beta=0$) to find for $r\in\Zp$
$$
\eqalign{
&\qquad\qquad \int_0^\infty z^{2r}
J_x(zq^{-{1\over 2}m};q)J_x(zq^{{1\over 2}k};q) dm_q(z) = \cr
&(-1)^{m+k}c^{{1\over {2n}}(m-k)}c^{-r} c^{-{k\over n}r}
\int_0^\infty z^{2r} p_{n+m+k}(z;n)p_n(z;n) dm_{c^{1\over n}}(z),\cr}
\tag5.2
$$
by shifting the summation parameter.
The orthogonality relations \thetag{3.1} imply
$$
\int_0^\infty p_k(z;n)p_l(z;n) dm_{c^{1\over n}}(z) = \delta_{k,l} .
$$
Repeating the first part of the proof of theorem 2 of Van Assche and
Koornwinder \cite{17} results in
$$
\lim_{n\to\infty} \int_0^\infty z^{2r} p_n(z;n)p_{k+n}(z;n)
dm_{c^{1\over n}} (z)
= {1\over{\pi}} \int_{B-2A}^{B+2A} {{z^rT_k\bigl( (z-B)/2A\bigr)}\over
{\sqrt{4A^2-(z-B)^2}}} \, dz
\tag5.3
$$
for $r,k\in\Zp$, $A=c^{1+\alpha}$, $B=c(1+c^{2\alpha})$ (as in \S 4) and
$T_k(\cos\theta)=\cos k\theta$ is the Chebyshev polynomial of the first kind
of degree $k$.

In \thetag{5.1} we replace $R$, $x$, $y$ and $q$ by
${1\over 2}R(1-q)$, $2n\alpha$, $2n\eta$, $c^{1\over n}$,
and we use \thetag{4.2}, \thetag{5.3},
\thetag{5.2} to see that \thetag{5.1}
(formally) tends to
$$
\eqalign{
&\qquad\qquad J_{\nu+m}(Rc^\eta) J_m(Rc^{\alpha+\eta}) = \cr
&\sum_{k=0}^\infty
{{(-\nu)_k}\over{\pi k!}} c^{\alpha k}(-1)^{k+m}
\int_{B-2A}^{B+2A} \bigl({z\over c}\bigr)^{-{1\over 2}\nu}
J_\nu(Rc^\eta \sqrt{{z\over c}})
{{T_{k+m}\bigl( (z-B)/2A\bigr)}\over{\sqrt{4A^2-(z-B)^2}}} dz\cr}
\tag5.4
$$
as $n\to\infty$. In the integral on the right hand side of
\thetag{5.4} we replace $(z-B)/2A$ by $\cos\psi$, so that
the right hand side of \thetag{5.4} equals
$$
\eqalign{
{1\over{2\pi}} \sum_{k=0}^\infty {{(-\nu)_k}\over{k!}} c^{\alpha k}
(-1)^{m+k}&\int_{-\pi}^{\pi}
(2c^\alpha\cos\psi+1+c^{2\alpha})^{-{1\over 2}\nu} \cr
& \times J_\nu\bigl( Rc^\eta \sqrt{2c^\alpha\cos\psi+1+c^{2\alpha}}\bigr)
e^{-i(k+m)\psi} \, d\psi.\cr}
\tag5.5
$$
Replace $\psi$ by $\psi+\pi$ in \thetag{5.5} and use
the binomial theorem to see that \thetag{5.5} equals
$$
{1\over{2\pi}} \int_0^{2\pi} J_\nu\bigl( Rc^\eta
\sqrt{1+c^{2\alpha}-2c^\alpha\cos\psi} \bigr)
{{(1-c^\alpha e^{-i\psi})^\nu}\over
{(1+c^{2\alpha}-2c^\alpha\cos\psi)^{{1\over 2}\nu}
}} e^{-im\psi}\, d\psi .
\tag5.6
$$
Equating the left hand side of \thetag{5.4} and
\thetag{5.6} yields a formula
equivalent to the product formula \thetag{1.6}.

\head Appendix. Justifications  \endhead

In the appendix we investigate the absolute convergence of two double
sums occuring in sections 2 and 3. The estimates used are based on the
estimates given in \cite{13}. In particular, cf. \cite{13, (2.4)},
$$
\vert J_\nu(z;q)\vert \leq \vert z^\nu\vert {{(-q^{\Re(\nu)+1},-q\vert
z^2\vert;q)_\infty}\over{(q;q)_\infty}}.
\tag{A.1}
$$
For integer order the Hahn-Exton $q$-Bessel function satisfies the
estimate
$$
\vert J_n(z;q)\vert \leq \vert z\vert^{\vert n\vert}
{{(-q,-q\vert z^2\vert;q)_\infty}\over{(q;q)_\infty}}
\cases 1, &\text{$n\geq 0$;}\\
       q^{{1\over 2}n(n-1)}, &\text{$n\leq 0$,}\endcases
\tag{A.2}
$$
by \cite{13, (2.4), (2.6)}. A similar estimate also holds for the
Hahn-Exton $q$-Bessel function with integer or half-integer power of $q$ as
argument because of the symmetry, cf. \cite{13, (2.3)},
$$
J_\alpha(q^{{1\over 2}\nu};q) = J_\nu(q^{{1\over 2}\alpha};q).
\tag{A.3}
$$

\proclaim{Proposition A.1}
The double sum
$$
\sum_{n=0}^\infty {{(-1)^nq^{{1\over 2}n(n+1)} R^{2n}q^{n(y+\nu)}}\over
{(q^{\nu+1};q)_n(q;q)_n}} \sum_{z=-\infty}^\infty q^{z(n+1+{1\over 2}\nu)}
J_x(q^{{1\over 2}(m+z)};q)J_{x-\nu}(q^{{1\over 2}z};q)
$$
is absolutely convergent for $m\in\Z$, $\Re (x) > -1$,
$\vert R\vert^2q^{1+\Re(x)+\Re(y)}<1$.
\endproclaim

\demo{Proof} Assume $m\geq 0$ and use \thetag{A.2}
and \thetag{A.3} to estimate
$$
\sum_{z=0}^\infty \vert q^{z(n+1+{1\over 2}\nu)}
J_x(q^{{1\over 2}(m+z)};q)J_{x-\nu}(q^{{1\over 2}z};q) \vert \leq
C \sum_{z=0}^\infty q^{z(1+n+\Re(x))} \leq C
$$
for $\Re(x)>-1$, $\forall\, n\in\Zp$ and $C$ is a constant independent
of $n$ of which the values may differ with each occurence.
Consequently, the part $\sum_{n=0}^\infty\sum_{z=0}^\infty$
is absolutely convergent for $\Re(x)>-1$.

Similarly we estimate
$$
\sum_{z=-\infty}^{-m} \vert q^{z(n+1+{1\over 2}\nu)}
J_x(q^{{1\over 2}(m+z)};q)J_{x-\nu}(q^{{1\over 2}z};q) \vert \leq
C \sum_{z=m}^\infty q^{z(1-n-m+\Re(x-\nu))}q^{z(z-1)}.
\tag{A.4}
$$
Shift the summation parameter to run from
$0$ to $\infty$, estimate $q^{z(z-1)}$
by $q^{{1\over 2}z(z-1)}/(q;q)_z$ and use the ${}_0\varphi_0$ summation
formula to see that \thetag{A.4} can be estimated by
$$
Cq^{-nm}(-q^{\Re(x-\nu)-n+1+m};q)_n \leq C q^{-{1\over 2}n(n-1)}
q^{n\Re(x-\nu)}.
$$
So the part $\sum_{n=0}^\infty \sum_{z=-\infty}^{-m}$ is absolutely
convergent for $\vert R\vert^2q^{1+\Re(x)+\Re(y)}<1$.

The remaining sum over $z$ is finite and causes no problems. The case
$m\leq 0$ is essentially the same. \qed
\enddemo

\proclaim{Proposition A.2}
The double sum
$$
\sum_{m=-\infty}^\infty q^m J_x(q^{{1\over 2}(l+m)};q)
\sum_{z=-\infty}^\infty q^z
J_x(q^{{1\over 2}(m+z)};q)J_{x-\nu}(q^{{1\over 2}z};q)
J_\nu(Rq^{{1\over 2}(y+\nu+z)};q)
$$
is absolutely convergent for $\Re (x) > -1$,
$\vert R\vert^2 q^{1+\Re(x)+\Re(y)}<1$.
\endproclaim

\demo{Proof} First we estimate the inner sum over $z$ as a function of $m$.
For the Hahn-Exton $q$-Bessel function of order $x$ and $x-\nu$ we use
the same estimates as before, and from \thetag{A.1} we get
$$
\vert J_\nu(Rq^{{1\over 2}(y+\nu+z)};q)\vert
\leq C q^{{1\over 2}\Re(\nu)z}
\cases 1,&\text{$z\geq 0$;} \\
\vert R\vert^{-2z} q^{-z(\Re(y+\nu)+1)}q^{-{1\over 2}z(z-1)},
&\text{$z\leq 0$,} \endcases
$$
where $C$ is independent of $m$. Denote by $S_m(1)$ the inner sum over
$z$ from $\max(0,-m)$ to $\infty$, then we find, using these estimates,
for $\Re(x) >-1$,
$$
\vert S_m(1)\vert \leq C \cases q^{{1\over 2}m\Re(x)}, &\text{$m\geq 0$;}\\
q^{-m-{1\over 2}m\Re(x)}, &\text{$m\leq 0$.}\endcases
$$
For $S_m(2)$, the inner sum over $z$ from $-\infty$ up to $\min(0,-m)$, we
obtain
$$
\vert S_m(2)\vert \leq C \cases q^{m({1\over 2}\Re(x)+\Re(y))}
\vert R\vert^{2m}, &\text{$m\geq 0$;} \\
q^{-{1\over 2}m\Re(x)+{1\over 2}m(m-1)}, &\text{$m\leq 0$.}\endcases
$$
The remaining finite sum over $z$ from $\min(0,-m)+1$ to
$\max(0,-m)-1$, denoted by $S_m(3)$, can be estimated by
$$
\vert S_m(3)\vert \leq C \cases q^{{1\over 2}m\Re(x)}
\bigl(1-(\vert R\vert^2q^{\Re(y)})^{m-1}\bigr), &\text{$m\geq 0$;}\\
q^{-{1\over 2}m\Re(x)}, &\text{$m\leq 0$.} \endcases
$$

To estimate the double sum of the proposition we combine
\thetag{A.2} and \thetag{A.3} to
obtain an estimate on $J_x(q^{{1\over 2}(l+m)};q)$. The sum
$m=-\infty$ to $\min(0,-l)$ and the finite sum from $m=\min(0,-l)+1$
to $\max(0,-l)-1$ are absolutely convergent.
The sum $m=\max(0,-l)$ to $\infty$ is absolutely convergent
for $\Re(x)>-1$ and
$q^{1+\Re(x)+\Re(y)}\vert R\vert^2<1$. \qed
\enddemo


\enddocument